%% file: main.tex
\documentclass[letterpaper,10pt,conference]{ieeeconf}

\IEEEoverridecommandlockouts
\usepackage{cite}
\usepackage{amsmath,amssymb,amsfonts}
\usepackage{algorithmic}
\usepackage{graphicx}
\usepackage{textcomp}
\usepackage{enumerate}

\usepackage{xcolor}
\usepackage{subcaption}
\usepackage{bm}
\usepackage{arydshln}
\usepackage[implicit=false, colorlinks=true, urlcolor=blue]{hyperref}

\newtheorem{theorem}{Theorem}

\newtheorem{proposition}{Proposition}

\newtheorem{example}{Example}

\newtheorem{cor}{Corollary}

\newcommand{\expectation}{\mathbb{E}}
\newcommand{\expectationp}[2][]{\mathbb{E} \ifx #1 \undefined \else _{#1} \fi \left[#2\right]}

\newcommand{\variancep}[2][]{\mathrm{Var} \ifx #1 \undefined \else _{#1} \fi \left[#2\right]}
\newcommand{\dataset}{\mathcal{D}}

\newcommand{\condbar}{\,|\,}  

\newcommand{\normaldist}{\mathcal{N}}

\DeclareMathOperator*{\argmin}{argmin} 
\DeclareMathOperator*{\argmax}{argmax} 

\title{\LARGE \bf On a closed-loop identification challenge in feedback optimization}
\author{Kristian Lindbäck Løvland$^{1, 2, *}$, Lars Struen Imsland$^{1}$ and Bjarne Grimstad$^{1, 2}$\thanks{*: kristian.lovland@ntnu.no}\thanks{1: Norwegian University of Science and Technology, Trondheim, Norway.}\thanks{2: Solution Seeker AS, Oslo, Norway.}
\thanks{© 2025 IEEE. Personal use of this material is permitted. Permission from IEEE must be obtained for all other uses, in any current or future media, including reprinting/republishing this material for advertising or promotional purposes, creating new collective works, for resale or redistribution to servers or lists, or reuse of any copyrighted component of this work in other works.}
}

\begin{document}
\maketitle
\thispagestyle{empty}
\pagestyle{empty}

\begin{abstract}
Feedback optimization has emerged as an effective strategy for steady-state optimization of dynamical systems. By exploiting models of the steady-state input-output sensitivity, methods of this type are often sample efficient, and their use of feedback ensures that they are robust against model error. Still, this robustness has its limitations, and the dependence on a model may hinder convergence in settings with high model error. We investigate here the effect of a particular type of model error: bias due to identifying the model from closed-loop data. Our main results are a sufficient convergence condition, and a converse divergence condition. The convergence condition requires a matrix which depends on the closed-loop sensitivity and a noise-to-signal ratio of the data generating system to be positive definite. The negative definiteness of the same matrix characterizes an extreme case where the bias due to closed-loop data results in divergence of model-based feedback optimization.
\end{abstract}

\section{Introduction}
\subsection{Feedback optimization}
Feedback optimization has gained popularity in recent years as a strategy for optimization of the steady-state operating points of dynamical systems \cite{hauswirth_optimization_2024}. As the name suggests, this class of methods works by treating setpoint selection as an optimization problem; this optimization problem is solved in an online fashion using continuous output measurements, providing the feedback mechanism suggested by the name.

In many settings, feedback optimization provides a good tradeoff between classical supervisory feedback control and offline model-based optimization. The formulation as an optimization problem means that setpoints can be defined implicitly in terms of an (e.g. economical) objective function, while the feedback which is introduced by optimizing around online output measurements ensures constraint satisfaction and resilience against unmeasured disturbances. Furthermore, feedback optimization typically provides less conservative behavior than offline methods based on robust or stochastic optimization.

Though model-free variants exist, feedback optimization algorithms typically rely on knowledge about the steady-state input-output sensitivity. If a high-fidelity model which provides this sensitivity is available, feedback optimization will typically converge under quite weak assumptions, and it will converge in fewer iterations than its model-free counterparts \cite{hauswirth_optimization_2024, he2024gray}.

Even in settings where such a model is not available, one will often have access to historical data from which a model can be learned. However, historical data must be applied with care. In settings where one wants to implement feedback optimization, an alternative (and presumably unsatisfactory) optimization strategy is likely to have been in use during data generation. If this strategy makes use of any type of feedback, estimation of the input-output sensitivity becomes a problem of \textit{closed-loop identification.} As is well known in system identification, statistical estimation of an input-output mapping performed directly on closed-loop data is not guaranteed to result in a model which represents the true causal input-output relation \cite{forssell_closed-loop_1999, van_den_hof_closed-loop_1998, ljung_system_1999}. Motivated by this, we analyze the performance of a class of feedback optimization methods when the steady-state input-output sensitivity is estimated from a type of closed-loop data.

Our main result characterizes conditions on the data generating system which ensure converge for a simple type of gradient-based feedback optimization method. We also present a converse condition under which the model bias is large enough for feedback optimization to diverge. The conditions apply to a steady-state LTI setting, and are stated in terms of properties of the closed-loop system from which data are observed, through a definiteness condition on a matrix defined from loop transfer functions and a signal-to-noise ratio. 

\subsection{Related work}
Feedback optimization is closely related to methods like simultaneous perturbation stochastic approximation \cite{spall1992multivariate, spall2000adaptive}, zeroth-order optimization \cite{liu2020primer}, Bayesian optimization \cite{shahriari2015taking, frazier2018tutorial}, extremum seeking \cite{ariyur2003real} and modifier adaptation \cite{chachuat2009adaptation, marchetti2009modifier, marchetti2016modifier}. In contrast to the formulation of feedback optimization considered in this article, however, these methods are typically not model-based. For a discussion of the relation between feedback optimization and all of these methods, as well as a broader discussion of feedback optimization methodology, see \cite{hauswirth_optimization_2024}. 

The effect of model error on feedback optimization was addressed in the context of bounded uncertainty sets by \cite{colombino_towards_2019}. Data-driven approaches to model-based feedback optimization have also been analyzed, e.g. by \cite{bianchin2023online}. Works like \cite{bianchin2023online}, which makes use of behavioral approach to system theory \cite{markovsky2021behavioral}, depend on open-loop input-output trajectories to build input-output sensitivity estimates. This is in contrast to our setup, where we only assume access to steady-state, closed-loop data.

The basic model alignment condition which is used in this paper bears resemblance with the model adequacy conditions which are often considered on the literature on real-time optimization and modifier adaptation \cite{Forbes1996}, as well as the more general condition of \textit{monotonicity}, which was applied to a problem of input-constrained feedback optimization by \cite{colombino_towards_2019}. While our results are derived on the context of online approximate gradient descent applied to unconstrained optimization, it is likely that the convergence conditions derived here can be extended to settings like these.

Finally, we note that our model convergence result is closely related to (and inspired by) standard results in closed-loop identification, see e.g. \cite{forssell_closed-loop_1999, van_den_hof_closed-loop_1998, ljung_system_1999}.



\subsection{Contributions}
This article investigates the robustness of feedback optimization to model error. There are many ways to mitigate this problem (for instance, one can always disregard the model completely and use a model-free method). Our investigation concerns a specific, particularly simple type of gradient-based feedback optimization strategy which has good convergence properties if the input-output model is good. Our results show how the fast convergence provided by a good input-output model can turn into divergence if the model is poor.

To our knowledge, the effect of closed-loop data on feedback optimization has not been addressed before. Our main contribution lies in relating a classical type of result from closed-loop system identification with a convergence condition for a simple but common gradient-based feedback optimization strategy.
The result is a description of conditions where a static, data-driven model of the steady-state input-output sensitivity is sufficient for feedback optimization to converge, as well as a description of particularly challenging conditions where this approach is too simple, leading to divergent behavior.

\subsection{Notation}
We denote the Jacobian of a function $f(x)$ as $\partial_x f(x)$, and its gradient as $\nabla_x f(x) = \partial_x f(x)^T$. For composed functions $f(x(z))$ evaluated at the point $(x', z')$, we abuse notation somewhat and write $\partial_x f(x') := \partial f(x)/\partial x \rvert_{x=x'}$ and $\partial_z f(x') = d f(x(z)) / dz\rvert_{z=z'}$.

\section{Convergence of feedback optimization}
\label{sec:convergence-and-monotonicity}
\subsection{General setup}
A typical problem to which feedback optimization is applied is the static optimization problem
\begin{equation}
    \min_{u \in \mathcal U} \phi(y) \quad \textup{s.t. } y = \pi(u) + w \label{eq:general-problem-statement}
\end{equation}
where $u \in \mathcal U \subseteq \mathbb R^{n_u}$ is the input, $y \in \mathbb R^{n_y}$ is the output, $w \in \mathbb R^{n_w}$ is an unobserved disturbance, $\phi : \mathbb R^{n_y} \rightarrow \mathbb R$ is the known cost function, and both $\phi$ and $w$ may be time varying. The map $\pi : \mathbb R^{n_u} \rightarrow \mathbb R^{n_y}$ is defined as the steady-state of a dynamical system, and its existence and uniqueness is usually proved through the use of an implicit function theorem. The treatment of $y$ as a static function of $u$ and $w$ is based on timescale separation arguments, commonly proved through the use of singular perturbation analysis; we refer to \cite{hauswirth2020timescale} for a discussion of interaction and separation between fast system dynamics and slow optimization dynamics. For simplicity, the cost function is often defined only in terms of the output $y$, but it is typically possible to add a term depending on the input $u$ without technical difficulties. For a broader description of problems like (\ref{eq:general-problem-statement}), see \cite{hauswirth_optimization_2024, krishnamoorthy2022real}.

\subsection{Linear unconstrained setup}
The most impactful (and arguably also the most enlightening) analyses of closed-loop data issues have appeared in the context of linear systems \cite{ljung_system_1999, van_den_hof_closed-loop_1998, forssell_closed-loop_1999}. Motivated by this, we consider a simplified version of (\ref{eq:general-problem-statement}) which lets us apply linear system theory.
We consider the static case, meaning $\phi$ and $w$ are not time-varying (but we note that results from the static case can often be applied to the case where $\phi$ and $w$ are time varying; see \cite{belgioioso2021sampled} for a discussion on how). We also assume that $n_y = n_u = n$, i.e. the system is square. Furthermore, we assume that there are no input constrains, meaning $\mathcal U = \mathbb R^n$, and we assume that $\phi$ is differentiable and strictly convex. Finally, we assume that $\pi(u) = \Pi u$ is linear, where $\Pi$ has full rank (this assumption is necessary since later derivations will rely on inverting $\Pi$). We can then restate (\ref{eq:general-problem-statement}) as:
\begin{equation}
    \min_{u \in \mathbb  R^n} \phi(y) \quad \textup{s.t. } y = \Pi u + w
\label{eq:linear-problem-statement}
\end{equation}
where $u, y \in \mathbb R^n$ are the input and output, $\phi : \mathbb R^n \rightarrow \mathbb R$ is the strictly convex cost function and the full rank matrix $\Pi \in \mathbb R^{n \times n}$ defines the true input-output mapping.

The problem formulation (\ref{eq:linear-problem-statement}) concerns unconstrained optimization. We do however note that the problem statement can also be relevant to settings of constrained optimization, for instance in settings where constraints can be eliminated through the use of active constraint control \cite{krishnamoorthy2022real}. Strict convexity of $\phi$ ensures convergence towards a global optimum; under the assumption that $\phi$ is differentiable but non-convex, the results derived later will remain valid, but convergence will only hold locally.

\subsection{Stability, convergence and monotonicity}
We investigate the convergence of a simple method for feedback optimization: The \textit{Online Approximate Gradient} (OAG) method (see e.g. \cite{colombino_towards_2019}). The OAG method can be defined in continuous-time as
\begin{equation}
    \dot u = - \varepsilon\hat g(u)
\end{equation}
where $\varepsilon > 0$ is a scaling factor and $\hat g(u)$ is a gradient estimate. Applied to (\ref{eq:linear-problem-statement}), this estimate would ideally equal the gradient of $\phi$ with respect to $u$:
\begin{align}
g(u) & = \nabla_u \phi(\pi(u) + w) \\
& = \partial_u \pi(u)^T \nabla_y \phi(y) \\
& = \Pi^T \nabla_y \phi(y)
\end{align}
In practice, however, the gradient is usually not known exactly. Since $\phi$ is known and $\Pi$ is unknown, a gradient estimate can be constructed from a model $\hat \Pi$ of the input-output sensitivity:
\begin{equation}
    \hat g(u) = \hat \Pi^T \nabla_y \phi(y)
    \label{eq:gradient-estimate}
\end{equation}
In practice one would typically implement OAG through its Euler approximation, i.e. as $u_{t+1} = u_t - \tau \hat g(u_t)$. For purposes of analysis, however, it is useful to consider the continuous-time formulation; the results will be applicable to the discrete-time formulation under suitable assumptions on the step size $\tau$ (see e.g. \cite{ljung_analysis_1977, metivier_applications_1984, borkar2008stochastic}).

The quality of the optimization scheme will depend on the quality of the gradient estimate. The following result provides a tool for analyzing how $\hat g(u)$ affects OAG convergence.

\begin{proposition}[Model alignment]
\label{prop:stability}
Let
\begin{align}
\begin{split}
    \dot u & = - \varepsilon \hat \Pi^T \nabla_y \phi(y) \\
    y & = \Pi u + w
    \label{eq:autonomous-system}
\end{split}
\end{align}
Assume that
$\phi$ is strictly convex and radially unbounded.
Suppose further that
\begin{equation}
    \Pi \hat \Pi^T + \hat \Pi \Pi^T \succ 0
    \label{eq:linear-monotonicity-condition}
\end{equation}
Then, (\ref{eq:autonomous-system}) is globally asymptotically stable, with equilibrium point $u = u^{\star} = \argmin_{u \in \mathbb R^n}\phi(\Pi u + w)$.
\end{proposition}
\begin{proof}
Since $\phi$ is assumed to be strictly convex (and, consequently, positive definite) and radially unbounded, it will be a Lyapunov function for the system (\ref{eq:autonomous-system}) if $\dot \varphi < 0$ for all $t$. Assuming without loss of generality that $\varepsilon=1$ and applying the chain rule, we get the condition
\begin{align}
    \frac{d}{dt} \phi(\Pi u + w) & = \nabla_u \phi(\Pi u + w)^T \dot u \\
    & = ( \Pi^T \nabla_y \phi(y) )^T (-\hat \Pi^T \nabla_y \phi(y)) \\
    & = - \nabla_y \phi(y)^T  \Pi \hat \Pi^T \nabla_y \phi(y) \\
    & \leq 0
\end{align}
If (\ref{eq:linear-monotonicity-condition}) holds, this will be satisfied for all $y \in \mathbb R^n$, with equality if and only if $\nabla_y \phi(y) = 0 \iff y = y^{\star} \iff u = u^{\star}$ (where the last equality follows from the fact that $\Pi$ has full rank). We conclude that (\ref{eq:linear-monotonicity-condition}) is a sufficient condition for the system (\ref{eq:autonomous-system}) to be globally asymptotically stable.
\end{proof}
Proposition \ref{prop:stability} illustrates the advantages of the OAG method: If one has access to the cost function $\phi(y)$ as well as a model estimate $\hat \Pi$ which satisfies the positive definiteness condition (\ref{eq:linear-monotonicity-condition}), the steady-state optimization problem (\ref{eq:linear-problem-statement}) can be solved with little additional effort. Due to the convex, unconstrained setup considered here, convergence towards the optimum is guaranteed even in the presence of significant model error.\footnote{The reason why this is possible is that we consider the setting of unconstrained inputs. In the case of input constraints, this result would be weaker; see e.g. \cite{colombino_towards_2019}, which showed convergence towards an approximate solution under monotonicity assumptions similar to (\ref{eq:linear-monotonicity-condition})}


The question we need to answer is then: \textit{When will a model $\hat \Pi$ learned from closed-loop data satisfy (\ref{eq:linear-monotonicity-condition})?}

\section{Steady-state closed-loop identification}
\label{sec:closed-loop-identification}
We proceed by characterizing the bias of $\hat \Pi$ in terms of closed-loop data generating system properties.

\subsection{System, controller and disturbances}
We assume that the system to be optimized is LTI and subject to additive noise. In subsequent discussions it will be convenient to describe the system with a frequency domain representation, such that:
\begin{equation}
    y = G(s) u + w \label{eq:system-dynamics}
\end{equation}
We also assume that during data generation, the system is governed by an LTI feedback controller:
\begin{equation}
    u = K(s) [r - y] \label{eq:controller-dynamics}
\end{equation}
Finally, we assume that $G(s)$ and $K(s)$ are proper, $G(s)K(s)$ is strictly proper, and that the closed-loop system (\ref{eq:system-dynamics})-(\ref{eq:controller-dynamics}) is stable.

Borrowing terms from hierarchical control design (see e.g. \cite{skogestad2005multivariable}, Chapter 10), we assume that $K(s)$ operates in the \textit{supervisory} or \textit{optimization} layer, while $G(s)$ operates on the lowest level in the control hierarchy. This means that for unstable plants, $G(s)$ will include the regulatory control loops which provide stabilization. In such cases, the value $u$ provided by $K(s)$ will typically contain reference values for the regulatory controllers.

\begin{example}[Optimizing PI controller]
\label{ex:pi-control}
Consider the feedback optimization method
\begin{equation}
    \dot u_t = -\varepsilon \hat \Pi^T \nabla_y \phi_t(y_t)
\end{equation}
where the time-varying cost function has the form
\begin{equation}
    \phi_t(y_t) = \frac{c_1}{2} \| r_t - y_t \|^2 + \frac{c_2}{2} \| \dot r_t - \dot y_t \| ^2
\end{equation}
As shown in Appendix \ref{sec:fo-is-lti-controller}, this optimization strategy defines a reference tracking PI controller:
\begin{equation}
    u_t = \underbrace{\hat \Pi^T \left(\frac{\varepsilon_1 + \varepsilon_2 s}{s} \right)}_{K(s)} (r_t - y_t)
\end{equation}
\end{example}
As Example \ref{ex:pi-control} illustrates, the set of LTI controllers $K(s)$ which we consider here contains a certain class of feedback optimization methods. We do however note that there are clearly many optimization strategies which can not be realized as linear controllers; just as the assumption that $\pi(u)$ is linear in $u$, this assumption is made for reasons of simplicity.

\subsection{Timescale separation}
We assume that the dynamics of $G(s)$ are significantly faster than those of the disturbance $w$ and reference $r$. We also assume that the controller $K(s)$ operates on the same timescale as the feedback optimization method by which it will be replaced. If this feedback optimization method is to counteract the disturbance and track the reference, its dynamics should be faster than those of the disturbance and reference. We will assume here that this is the case.

Denoting the time constants of $G(s)$, $K(s)$, $r$ and $w$ somewhat qualitatively as $\tau_G$, $\tau_K$, $\tau_r$ and $\tau_w$, these assumptions can be stated compactly as $\tau_G, \tau_K \ll \tau_r, \tau_w$. The assumption $\tau_G \ll \tau_r, \tau_w$ is standard in feedback optimization. Since the optimization layer $K(s)$ treats the system layer $G(s)$ as a steady-state mapping, we will also have $\tau_G \ll \tau_K$. This assumption is compatible with, but not necessary for, the results which we will derive. The assumption $\tau_K \ll \tau_r, \tau_w$ can be motivated from works like \cite{colombino2019online}, which require that the feedback optimization dynamics are considerably faster than the disturbance and reference processes.

We stress that while the time constant of the existing controller $K(s)$ may be sufficiently fast to track $r_t$, we do not assume that this tracking is actually achieved. Imperfect control provided by $K(s)$ is a central motivation for its replacement by a feedback-based optimization method.


\subsection{Dataset}
We assume that we have access to a dataset $\dataset_T = \{ (u_t, y_t)\}_{t=1}^T$ of input-output observations arising from (\ref{eq:system-dynamics})-(\ref{eq:controller-dynamics}). By the timescale separation assumption $\tau_G, \tau_K \ll \tau_r, \tau_w$, these observations will be at steady-state:
\begin{align}
    y & = G(0) u + w \label{eq:system-steady-state} \\
    u & = K(0) [r - y] \label{eq:controller-steady-state}
\end{align}
From the timescale separation assumption it also follows that $r_t$ and $w_t$ can be viewed as constants for each observation time $t \in \{1, \dots, T\}$.

We emphasize that $\Pi = G(0)$, and that (\ref{eq:system-steady-state}) and the constraint of (\ref{eq:linear-problem-statement}) represent the same system. In the following, we will use $\Pi$ when discussing feedback optimization, and we will denote $G(0)$ compactly as $G$ when discussing the data generating system. Similarly, we will denote $K(0)$ as $K$. We assume that $G$ and $K$ are invertible, and that $G \neq - K^{-1}$.

We also emphasize that $\dataset_T$ consists of closed-loop, \textit{historical} data. Consequently, it grows backwards in time when $T$ is increased, and the limit case $T \rightarrow \infty$ represents an asymptotic setting where one has access to an infinitely long historical record of steady-states as described by (\ref{eq:system-steady-state})-(\ref{eq:controller-steady-state}).

We note that the assumptions which are made here concerning availability of offline data are weaker than what is common in the literature on feedback optimization; we assume here that only steady-state data are available, while datasets which are used for offline model estimation in the literature on feedback optimization is typically assumed to contain transient data, as is e.g. done by \cite{bianchin2023online}.


We assume that the input-output steady-state sensitivity model $\hat y = \hat \Pi_T u$ is identified from $\dataset_T$ through direct MSE minimization, or equivalently, that it is identified through maximum likelihood estimation of the probabilistic model $q_T(y \condbar u) = \normaldist(\hat \Pi_T u, \sigma^2)$ where $\sigma > 0$ is fixed. Thus, we assume that the disturbance $w_t$ is not modeled. This assumption greatly simplifies modeling (and is therefore rather common), but we will see that it is closely related to the risk for feedback optimization divergence.
 
We assume that the sequences $\{ w_t \}$ and $\{ r_t \}$ are uncorrelated, and that they are ergodic stationary with stationary distributions given by zero-mean, Gaussian distributions with covariance $\Sigma_w$ and $\Sigma_r$, respectively:
\begin{align}
    r & \sim \normaldist(0, \Sigma_r) \label{eq:reference-stationary-distribution} \\
    w & \sim \normaldist(0, \Sigma_w) \label{eq:disturbance-stationary-distribution}
\end{align}
This property is e.g. satisfied by stable colored noise sequences driven by zero-mean Gaussian white noise.

\subsection{Model convergence}
We now consider the limiting behavior of the model $\hat \Pi_T$ as $T \rightarrow \infty$. This setting represents a best-case scenario with regards to data availability: Since one has access to an infinitely long historical record of data, model error can be solely attributed to bias.

Defining the \textit{closed-loop sensitivity} $S = (I + G K)^{-1}$, it follows from (\ref{eq:system-steady-state})-(\ref{eq:controller-steady-state}) that the steady-state of the closed-loop system satisfies
\begin{equation}
    \begin{bmatrix}
        y \\ u
    \end{bmatrix}
    =
    \begin{bmatrix}
        S & S G K \\
        - K S & K S
    \end{bmatrix}
    \begin{bmatrix}
        w \\ r
    \end{bmatrix}
    \label{eq:exogenous-to-endogenous-mapping}
\end{equation}
Since $(u, y)$ is a linear function of the ergodic stationary zero-mean Gaussian random variable $(r, w)$, the random variable $(u, y)$ is also ergodic stationary, zero-mean and Gaussian. Furthermore, the distribution $p(u, y)$ from which our dataset is drawn can be written in closed form as a Gaussian distribution. Define now
\begin{equation}
\Lambda := \Sigma_r [\Sigma_r + \Sigma_w]^{-1} \label{eq:signal-to-noise-ratio-def}
\end{equation}
For a signal $r$ subject to additive noise $w$, the matrix $\Lambda$ provides the linear least-squares (or \textit{Wiener} filtered) reconstruction $\hat r = \Lambda d$ of the signal $r$ given the observation $d = r + w$ \cite{kamen2012introduction}. In the scalar case, $\Lambda = \frac{\sigma_r^2}{\sigma_r^2 + \sigma_w^2}$ is a normalized ratio between the reference variance and the variance of the disturbance $d = r + w$.  Motivated by this, we will call $\Lambda$ a \textit{signal-to-noise} ratio. 

Having defined $\Lambda$, we can derive a characterization of the limiting behavior of $\hat \Pi_T$.

\begin{theorem}[Model convergence]
\label{thm:model-convergence}
Assume that $\{ r_t \}_{t=1}^T$ and $\{ w \}_{t=1}^T$ are ergodic stationary with stationary distributions given by (\ref{eq:reference-stationary-distribution})-(\ref{eq:disturbance-stationary-distribution}), where $\Sigma_r, \Sigma_w \succ 0$. Furthermore, assume that $(u_t, y_t)$ satisfies (\ref{eq:system-steady-state})-(\ref{eq:controller-steady-state}) for all $t \in [1, T]$. Then, the ML estimator
\begin{equation}
    \hat \Pi_T = \argmax_{\Pi \in \mathbb R^{n \times n}} \frac{1}{T} \sum_{t=1}^T \log \normaldist(y_t; \Pi u_t, \sigma^2)
    \label{eq:model-estimation}
\end{equation}
where $\sigma > 0$ is fixed, satisfies
\begin{equation}
\hat \Pi_T \xrightarrow{T\rightarrow\infty} \Lambda \Pi + (I - \Lambda) (- K^{-1}) \quad \textup{a.s.}
\end{equation}
\end{theorem}
\begin{proof}
It follows from standard ML consistency results for linear-Gaussian models (see e.g. Example 7.8 of \cite{hayashi2011econometrics}) that if $\expectation[u u^T]  \succ 0$, then $\hat \Pi_T u$ converges to the conditional expectation $ \expectation[y \condbar u] $ defined by the stationary distribution $p(u, y)$. As shown in Appendix \ref{sec:closed-form-gaussian},
\begin{equation}
    p
    \left(
    \begin{bmatrix}
        y \\ u
    \end{bmatrix}
    \right)
    =
    \normaldist
    \left(
    \begin{bmatrix}
        0 \\ 0
    \end{bmatrix}
    ,
    \begin{bmatrix}
        \Sigma_y & \Sigma_{uy} \\
        \Sigma_{uy}^T & \Sigma_u
    \end{bmatrix}
    \right)
    \label{eq:multivariate-gaussian-formula}
\end{equation}
where
\begin{align}
    \Sigma_u & = KS\Sigma_w (K S)^T + KS\Sigma_r (K S)^T \label{eq:sigma_u}\\
    \Sigma_{uy} & = -S\Sigma_w (K S)^T + SGK\Sigma_r (K S)^T \label{eq:sigma_uy}
\end{align}
The conditional distributions of multivariate Gaussians can be written in closed form (see e.g. \cite{petersen_matrix_2008}). For zero-mean variables the expression is
\begin{equation}
    p(y \condbar u) = \normaldist(\Sigma_{uy} \Sigma_u^{-1} u, \Sigma_y - \Sigma_{uy} \Sigma_u^{-1} \Sigma_{uy}^T)
\end{equation}
Since the full rank of $G$ and $K$ and positive definiteness of $\Sigma_r$ and $\Sigma_w$ implies that $\expectation[u u^T] = \Sigma_u \succ 0$, this implies that $\hat \Pi_T \rightarrow \Sigma_{uy} \Sigma_u^{-1}$. As shown in Appendix \ref{sec:conditional-gaussian},
\begin{align}
    \Sigma_{uy} \Sigma_u^{-1} = & \Sigma_r (\Sigma_w + \Sigma_r)^{-1} G \nonumber \\
    & + (I - \Sigma_r (\Sigma_w + \Sigma_r)^{-1}) (-K^{-1})
\end{align}
Insertion of (\ref{eq:signal-to-noise-ratio-def}) along with $G = \Pi$ yields the desired result.
\end{proof}
Theorem \ref{thm:model-convergence} mirrors classical results in closed-loop identification which describe model estimators that are unable to separate input-output correlations arising from open-loop dynamics $G(s)$ and correlations introduced by the feedback controller $K(s)$.
While there are some notable technical differences between the dynamical and steady-state setting, the problem is, in both cases, related to the fact that
the model estimator is derived under the assumption that the observed data were generated in open-loop. This assumptions is incorrect in the case of closed-loop data, since the observed correlation between $u$ and $y$ in the available dataset in this case can be partly attributed to the unobserved disturbance $w$. In statistical terms, $w$ acts as a \textit{confounder}. The failure to account for $w$ is the reason for the bias in $\hat \Pi_T$.

In settings with closed-loop data, consistent estimation of $\Pi$ may require significant modeling efforts. When the characteristics of the disturbance $w$ is unknown, consistent estimation may lie outside the realm of feasibility. We proceed by deriving convergence conditions which may aid model design in such cases.

\section{Stability conditions}
We now consider the behavior of OAG when using the model $\hat \Pi: = \lim_{T \rightarrow \infty} \hat \Pi_T$.
Having proved Theorem \ref{thm:model-convergence}, we can apply Proposition \ref{prop:stability} and derive a condition on the closed-loop system which ensures convergence of OAG which uses a model $\hat \Pi$ which is subject to closed-loop bias.
\begin{theorem}[OAG convergence]
\label{thm:oag-stability}
Assume that $\hat \Pi_T$ is estimated through MSE minimization, as described by (\ref{eq:model-estimation}), where the data arise from the closed-loop steady-state system (\ref{eq:system-steady-state})-(\ref{eq:controller-steady-state}). Suppose that this closed-loop system satisfies
\begin{equation}
    \frac{1}{2}(I - \Lambda)(K S)^{-1}G^T + \frac{1}{2}G (K S)^{-T}(I - \Lambda) \prec G G^T
    \tag{C}
    \label{eq:oag-stability-condition}
\end{equation}
where
$S = (I + G K)^{-1}$ and $\Lambda = \Sigma_r (\Sigma_r + \Sigma_w)^{-1}$. Assume that $w \in \mathbb R^n$ is constant and that $\phi$ is strictly convex and radially unbounded, and define $\hat \Pi = \lim_{T \rightarrow \infty} \hat \Pi_T$. Then, the system
\begin{align}
\begin{split}
    \dot u & = - \varepsilon \hat \Pi^T \nabla_y \phi(y) \\
    y & = \Pi u + w
    \label{eq:autonomous-system-restated}
\end{split}
\end{align}
is globally asymptotically stable, and its equilibrium point $u^{\star}$ satisfies $\Pi u^{\star} + w = \argmin_{y \in \mathbb R^n} \phi(y)$.
\end{theorem}
\begin{proof}
Since $S = (I + G K)^{-1}$ and $\Pi = G$, the limiting expression of Theorem \ref{thm:model-convergence} can be rewritten as
\begin{align}
    \hat \Pi & = \Lambda G + (I - \Lambda) (-K^{-1}) \\
    & = G - (I - \Lambda) (K^{-1} + G) \\
    & = G - (I - \Lambda)(I + GK)K^{-1} \\
    & = G - (I - \Lambda)(K S)^{-1} \label{eq:additive-bias-decomposition}
\end{align}
That is, $\hat \Pi = \Pi - B$, where $B = (I - \Lambda)(K S)^{-1}$. Under this decomposition, the condition $\Pi \hat \Pi^T + \hat \Pi \Pi^T \succ 0 $ is equivalent to $ \frac{1}{2}B G^T + \frac{1}{2}G B^T \prec G G^T$. Insertion of $B = (I - \Lambda)(K S)^{-1}$ and application of Proposition \ref{prop:stability} then yields the desired result.
\end{proof}
\begin{cor}[OAG divergence]
\label{cor:oag-instability}
Suppose that
\begin{equation}
    \frac{1}{2}(I - \Lambda)(K S)^{-1}G^T + \frac{1}{2}G (K S)^{-T}(I - \Lambda) \succ G G^T
    \tag{C'}
    \label{eq:oag-instability-condition}
\end{equation}
Then, the system (\ref{eq:autonomous-system-restated}) is unstable.
\end{cor}

By Theorem \ref{thm:oag-stability}, the definiteness condition (\ref{eq:oag-stability-condition}) ensures convergence of the continuous-time formulation of OAG. Under suitable assumptions on the step length (which must be decreasing at an appropriate rate), the same condition can be applied to discrete-time OAG formulations as well. The condition (\ref{eq:oag-instability-condition}) (which is not the logical negation of (\ref{eq:oag-stability-condition}), but a converse negative definiteness condition) describes an extreme case where the bias due to closed-loop data renders the OAG algorithm divergent.

\begin{figure}[bt]
    \centering
    \includegraphics[width=0.8\linewidth]{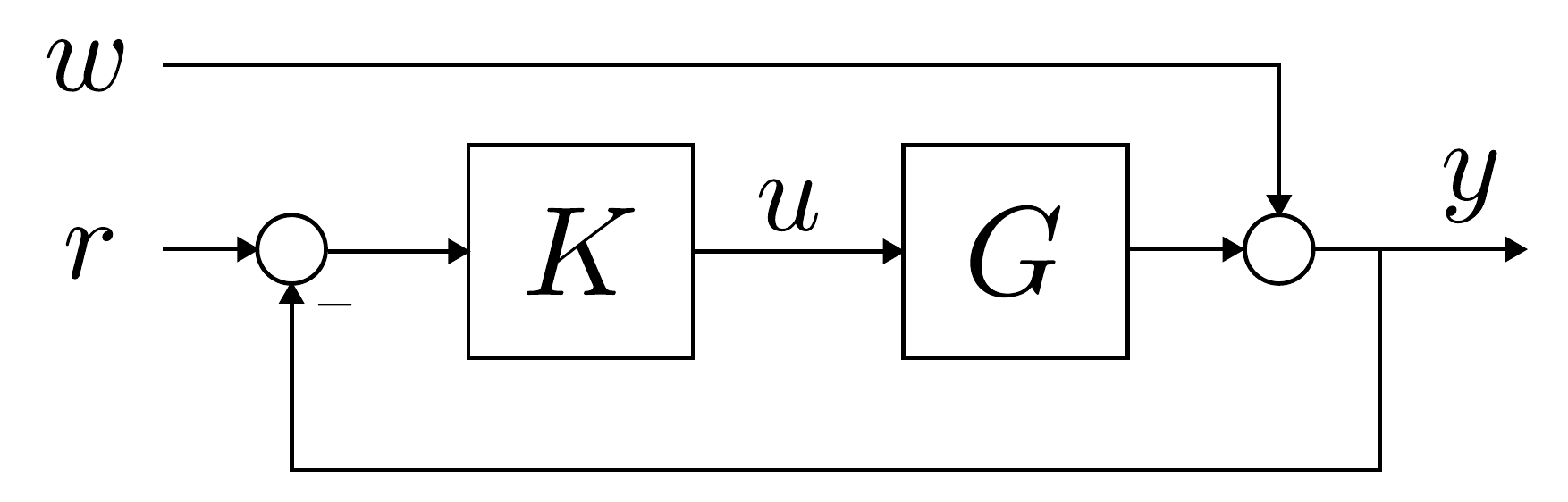}
    \caption{Structure of linear feedback controlled system}
    \label{fig:block-diagram}
\end{figure}

The additive bias term $B = (I - \Lambda)(K S)^{-1}$ lends itself to a fairly direct control theoretical interpretation. The mapping $K S$ represents the steady-state of the closed-loop transfer function from the reference $r$ to the closed-loop control input $u$. As can be seen from Figure \ref{fig:block-diagram}, the entry points of the exogenous terms $r$ and $-w$ can be merged into a single addition. Defining $d = r - w$, we can see from Figure \ref{fig:block-diagram} (and equivalently, from (\ref{eq:exogenous-to-endogenous-mapping})) that $u = K S d$. Since $K S$ is invertible, the difference $d = r - w$ between the exogenous terms can be calculated exactly from $u$ during closed-loop operation as $d = (K S)^{-1} u$. The matrix $(I - \Lambda)$ represents a static Wiener filter which provides a least-squares estimate $\hat w = (I - \Lambda) d$ of $w$ given $d$. The term $Bu$ can, then, be interpreted as a least-squares estimate of $w$ given a closed-loop observation of $u$.

The condition (\ref{eq:oag-stability-condition}) concerns model gradients.
If we let $\Delta y \in \mathbb R^n$ denote an arbitrary output perturbation (which in the context of OAG represents the desired output direction $\Delta y = \nabla_y \phi(y)$), and denote the least-squares disturbance estimate as $\hat w = B u$ and true input-output sensitivity as $y = \Pi u$, we can see from the definition of positive definite matrices that (\ref{eq:oag-stability-condition}) is equivalent with the condition that $\nabla_u \hat w$ is ``smaller'' than  $\nabla_u y$, in the sense that
\begin{equation}
    \langle \nabla_u \hat w, \nabla_u y \rangle < \langle \nabla_u y, \nabla_u y \rangle
\end{equation}
In other words, the information contained by $u$ about $w$ in closed-loop must be weaker than its effect on $y$. This must be the case for all output directions $\Delta y \in \mathbb R^n$.

Simplified conditions on the closed-loop data generating system can be derived by considering special cases of $\Lambda$ and $S$. If the reference $r$ is significantly stronger than the disturbance $w$ we will have $\Lambda \approx I$, and (\ref{eq:oag-instability-condition}) will simplify to $G G^T \succ 0$, which will always hold since $G$ is assumed to have full rank.
For a controller which achieves perfect reference tracking, the closed-loop steady-state will be completely insensitive to the disturbance (i.e. $S = 0$). In this case, one can show that (\ref{eq:oag-stability-condition}) holds as long as $\Lambda \succ 0$. \footnote{This follows from the fact that for a controller which approaches perfect reference tracking $y = r$ due to integral action in the controller, we will have $K^{-1} \rightarrow 0$. Since $(K S)^{-1} = (I + GK)K^{-1} = G + K^{-1}$, (\ref{eq:oag-stability-condition}) reduces to $\frac{1}{2}(I - \Lambda)GG^T + \frac{1}{2}G G^T(I - \Lambda) \prec G G^T$, which holds as long as $\Lambda \succ 0$.} This, in turn, is true if $\Sigma_r \succ 0$.

In the scalar case, i.e. $G, K, \Lambda \in \mathbb R$, (\ref{eq:oag-stability-condition}) can, under the assumption $K S > 0$, be simplified to
\begin{align}
    \frac{1 - \Lambda}{K S} & < G \\
    \iff 1 - \Lambda & < GKS = 1 - S \\
    \iff \Lambda & > S
\end{align}
In this case, the interpretation can be given directly in terms of $\Lambda$ and $S$: The static Wiener filter gain from $d$ to $r$ in closed-loop must be larger than the closed-loop sensitivity from $w$ to $y$.

\begin{figure}[bt]
    \centering
    \begin{subfigure}[b]{0.48\linewidth}
        \centering
        \includegraphics[width=\linewidth, trim=0.6cm 0cm 0.52cm 0cm, clip]{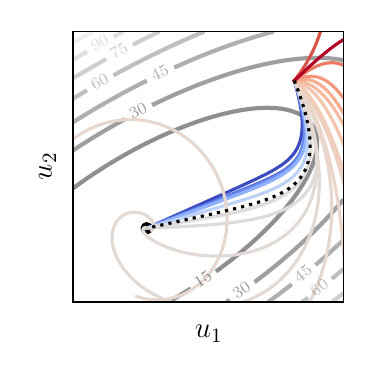}
        \caption{Control inputs}
        \label{fig:control-inputs}
    \end{subfigure}
    \begin{subfigure}[b]{0.50\linewidth}
        \centering
        \includegraphics[width=\linewidth, trim=0.6cm 0cm 0.3cm 0cm, clip]{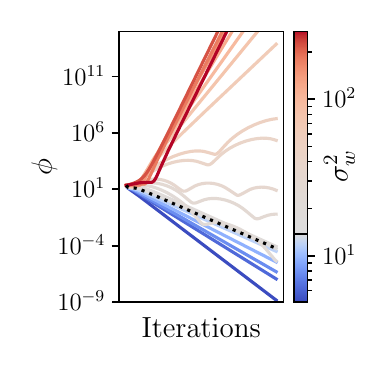}
        \caption{Cost function}
        \label{fig:cost-function}
    \end{subfigure}
    \caption{OAG in the setting described in Example \ref{ex:numerical-example}. The gray ellipsoids represent level curves of $\phi$, the colored lines in Figure \ref{fig:control-inputs} represent the trajectory of $u$, and the colored lines in Figure \ref{fig:cost-function} represent the corresponding value of the cost function $\phi$. The colors represent different noise levels $\sigma_w^2$, with higher colors corresponding to more noise. The black, dotted line in the figures and the black tick in the colorbar corresponds to the value of $\sigma_w^2$ for which (\ref{eq:oag-stability-condition}) holds with the smallest margin.}
    \label{fig:example}
\end{figure}

\begin{example}
\label{ex:numerical-example}
Consider a two-dimensional instance of (\ref{eq:linear-problem-statement}) where $\phi(y_1, y_2) = y_1^2 + y_2^2$ and
\begin{align}
\begin{split}
    \Pi & =
    \begin{bmatrix}
    1 & 2 \\
    -3 & 4
    \end{bmatrix},
    \quad
    K = 
    \begin{bmatrix}
    10 & 1 \\
    3 & 2
    \end{bmatrix}
    \\
    \Sigma_r & = 
    \begin{bmatrix}
    2 & 1 \\
    1 & 3
    \end{bmatrix},
    \quad
    \Sigma_w
    = 
    \begin{bmatrix}
    \sigma_w^2 & 0 \\
    0 & \sigma_w^2
    \end{bmatrix}
\end{split}
\end{align}
Figure \ref{fig:example} shows the result of applying OAG with models $\hat \Pi$ fit with least squares to 5000 i.i.d. closed-loop data points from these systems, starting in the point $u_0 = (-0.75, 1.5)$ and using a step size of $10^{-3}$. The condition (\ref{eq:oag-stability-condition}) holds for sufficiently low values of $\sigma_w^2$, with the threshold value being approximately 13.7. For all of the blue lines in the figure, (\ref{eq:oag-stability-condition}) holds, and the OAG algorithm converges. Some of the red lines (corresponding to lower values of $\sigma_w^2$) also converge towards the optimum, despite not satisfying (\ref{eq:oag-stability-condition}). This is expected, since (\ref{eq:oag-stability-condition}) is a sufficient, but not necessary, condition for convergence.
\end{example}

\section{Concluding remarks}
Using ideas from closed-loop system identification, we have investigated problems which arise when using historical data for model-based feedback optimization. Our main result can be used to analyze whether direct input-output sensitivity estimation from historical data will result in convergent feedback optimization, or if a less biased (and likely more involved) estimation scheme must be applied. The key condition which underpins this convergence result is the matrix definiteness condition (\ref{eq:oag-stability-condition}), which concerns properties of the closed-loop system from which data are observed. This condition holds if the gain from reference to output is, in a certain sense, larger than the fraction of the total input signal which arises from the disturbance. The converse definiteness condition (\ref{eq:oag-instability-condition}) presents an extreme case where model bias arising from closed-loop data leads to divergence of feedback optimization. The divergence condition illustrates that for overly simple input-output sensitivity estimates, closed-loop data issues may result in feedback optimization divergence, even in the presence of unlimited historical data.

We hope that the result can aid analysis and discussion of data issues when dealing with static optimization problems, offering some of the conceptual clarity which linear system theory often provides.

Concerning future work, a natural research direction would be to generalize the result to a less restrictive (e.g. nonlinear) setting. Another interesting problem to address would be a setting where the model $\hat \Pi_T$ is continually updated while being used for feedback optimization, and in particular, what excitation conditions would need to be enforced for such an adaptive optimization approach to work.

\bibliographystyle{IEEEtran}
\bibliography{bibliography} 
\appendix
\input{appendix}

\end{document}

%% file: appendix.tex
\section{Derivations}

\subsection{Derivation of LTI form of feedback optimization}
\label{sec:fo-is-lti-controller}
The gradient of $\phi_t(y_t)$ is given by
\begin{equation}
    \nabla_y \phi(y_t) = - c_1 (r_t - y_t) - c_2 (\dot r_t - \dot y_t)
\end{equation}
Let $\varepsilon_1 = \varepsilon c_1$, $\varepsilon_2 = \varepsilon c_2$. In the Laplace domain, the controller $u$ then satisfies
\begin{align}
    s u & = \Pi_0^T \left( \varepsilon_1 (r_t - y_t) + \varepsilon_2 (s r_t - s y_t) \right) \\
    \iff u & = \Pi_0^T \left( \frac{\varepsilon_1 + \varepsilon_2 s}{s} \right) (r_t - y_t)
\end{align}

\subsection{Derivation of $p(u, y)$}
\label{sec:closed-form-gaussian}
The linear transformation $z = A x$ of a zero-mean Gaussian variable $x \sim \normaldist(0, \Sigma)$ has distribution $p(z) = \normaldist(0, A \Sigma A^T)$. Insertion of (\ref{eq:exogenous-to-endogenous-mapping}) yields
\begin{align}
    & \begin{bmatrix}
        \Sigma_y & \Sigma_{uy} \\
        \Sigma_{uy}^T & \Sigma_u
    \end{bmatrix}
    \\
    &=
    \begin{bmatrix}
        S & SGK \\
        -KS & KS
    \end{bmatrix}
    \begin{bmatrix}
        \Sigma_w & 0 \\
        0 & \Sigma_r
    \end{bmatrix}
    \begin{bmatrix}
        S^T & - (K S)^T \\
        (S G K)^T & (K S)^T
    \end{bmatrix} \nonumber \\
    &=
    \begin{bmatrix}
        S\Sigma_w & SGK\Sigma_r \\
        -KS\Sigma_w & KS\Sigma_r
    \end{bmatrix}
    \begin{bmatrix}
        S^T & - (K S)^T \\
        (S G K)^T & (K S)^T
    \end{bmatrix}
\end{align}
The last matrix multiplication gives the desired result:
\begin{align}
    \Sigma_u & = KS\Sigma_w (K S)^T + KS\Sigma_r (K S)^T \\
    \Sigma_{uy} & = -S\Sigma_w (K S)^T + SGK\Sigma_r (K S)^T
\end{align}

\subsection{Derivation of $\Sigma_{uy} \Sigma_u^{-1}$}
\label{sec:conditional-gaussian}
Let $\Sigma = \Sigma_r + \Sigma_w$. Inserting the expressions from (\ref{eq:sigma_u}) and (\ref{eq:sigma_uy}) and performing some tedious but straight-forward matrix multiplication, we get the desired result:
\begin{align}
    \Sigma_{uy} \Sigma_u^{-1} & = \left[-S\Sigma_w (K S)^T + SGK\Sigma_r (K S)^T \right] \nonumber \\
    & \quad\quad\quad\left[ KS\Sigma_w (K S)^T + KS\Sigma_r (K S)^T \right]^{-1} \nonumber \\
    & = S [- \Sigma_w + G K \Sigma_r](K S)^T \nonumber \\ 
    & \quad\quad\quad\quad (K S)^{-T} [\Sigma_w + \Sigma_r]^{-1}(K S)^{-1} \nonumber \\
    & = S [- \Sigma_w + G K \Sigma_r]\Sigma^{-1}(K S)^{-1} \nonumber \\
    & = S [-\Sigma_w - \Sigma_r + \Sigma_r + G K \Sigma_r]\Sigma^{-1}(K S)^{-1} \nonumber \\
    & = S [- \Sigma + S^{-1} \Sigma_r]\Sigma^{-1} (K S)^{-1} \nonumber \\
    & = - K^{-1} + \Sigma_r \Sigma^{-1} S^{-1} K^{-1} \nonumber \\
    & = - K^{-1} + \Sigma_r \Sigma^{-1} (K^{-1} + G) \nonumber \\
    & = \Sigma_r \Sigma^{-1} G + (I - \Sigma_r \Sigma^{-1}) (-K^{-1}) 
\end{align}

